\documentclass[runningheads]{llncs}
\usepackage{graphicx}
\usepackage{algorithm}
\usepackage{algpseudocode}
\usepackage{subcaption}
\usepackage{a4wide}

\usepackage{amsmath}
\usepackage{amssymb}
\usepackage{mathtools}

\usepackage{booktabs}

\usepackage{adjustbox}

\usepackage{multirow}
\usepackage{tabularx}

\newcommand{\rank}{\mathrm{rank}}

\def\R{\mathbb{R}}
\def\RR{\mathbb{R}}

\newcommand{\norm}[1]{\left\|#1\right\|}

\def\argmax{\mathop{\rm arg\,max}\limits}%    a math operator.
\def\argmin{\mathop{\rm arg\,min}\limits}%    a math operator.

\def\min{\mathop{\rm min}\nolimits}
\def\max{\mathop{\rm max}\nolimits}

\begin{document}
\title{On Riemannian Approach for Constrained Optimization Model in Extreme Classification Problems\thanks{Supported by  IIIT, Hyderabad}
}
%
%\titlerunning{Abbreviated paper title}
% If the paper title is too long for the running head, you can set
% an abbreviated paper title here
%
\author{Jayadev Naram\inst{1} \and
Tanmay Kumar Sinha\inst{2} \and
Pawan Kumar\inst{3}}
%
%\authorrunning{Naram, J. et. al.}
% First names are abbreviated in the running head.
% If there are more than two authors, 'et al.' is used.
%
\institute{Internation Institute of Information Technology, Hyderabad, 500032, India \\
\inst{1}\email{jayadev.naram@students.iiit.ac.in} \\
\inst{2}\email{tanmay.kumar@research.iiit.ac.in}\\
\inst{3}\email{pawan.kumar@iiit.ac.in}}

\maketitle              % typeset the header of the contribution
\begin{abstract}
% A Riemannian solver for the global and sparse local embedding models for extreme classification problem is studied. The proposed solver improves on the state-of-the-art leading to a significant reduction in train time, while not exceeding the model size.  
We propose a novel Riemannian method for solving the Extreme multi-label classification problem that exploits the geometric structure of the sparse low-dimensional local embedding models. A constrained optimization problem is formulated as an optimization problem on matrix manifold and solved using a Riemannian optimization method. The proposed approach is tested on several real world large scale multi-label datasets and its usefulness is demonstrated through numerical experiments. The numerical experiments suggest that the proposed method is fastest to train and has least model size among the embedding-based methods. An outline of the proof of convergence for the proposed Riemmannian optimization method is also stated. 


\end{abstract}
\section{Introduction}
A generalization of multi-class classification is multi-label classification, where each data sample can belong to more than one class. Essentially, if we have $l$ labels, then each sample can be classified into as a subset of these, leading to $2^l$ distinct possibilities. Extreme multi-label classification problems focus on multi-label classification problems that involve a very large number of labels, sometimes more than a million. It is an important research problem with numerous applications, in tagging, ranking and recommendation. 

A major challenge faced in extreme classification problems is the problem of dealing with the tail labels, that is the labels that appear rarely in the dataset. To deal with this challenge, the {\tt SLEEC} \cite{sleec_bhatia} algorithm was proposed, that uses low-dimensional local embeddings to deal with the tail labels. The optimization model proposed by {\tt SLEEC} has a rich geometric structure; in this paper, we exploit this structure to develop a Riemannian variant of {\tt SLEEC}, which utilizes optimization techniques on matrix manifolds. This has a much faster training time without losing out on the predictive performance, and scales well with large datasets. We compare our algorithm with many proposed extreme classification algorithms, and detail the results. The main contributions of the paper are listed below:

% For bullet points in the list
\let\labelitemi\labelitemii
\begin{itemize}
    \item We propose a novel embedding-based Riemannian algorithm for extreme classification problems that exploits the rich geometric structure of the optimization model. We also outline a proof of convergence for the algorithm.
    \item The algorithm has significantly faster training time when compared to many state-of-the-art embedding methods, and fairly comparable train times with Tree-based methods. It also has a comparable predictive performance and the model size is least among all the embedding methods.
\end{itemize}

The rest of this paper is organized as follows. In section \ref{notation}, we introduce the notation used in the paper, and in section \ref{previous_work}, we review the previous work in the field of Extreme multi-label learning. In Section \ref{optimization_model}, we discuss the optimization model in detail, and present the basics ideas of optimization on matrix manifolds. We then present the Riemannian algorithm in full detail, discuss the space and time complexity, and also give a tentative proof of the convergence of the algorithm. In Section \ref{numerical_experiments} we detail the experiments we carried out and compare our algorithm to other well known algorithms.
In Section 6 we end with concluding remarks.

%\section{Problem Formulation}
% Extreme multi-label learning is an important research problem and it has several applications
% in tagging, recommendation and ranking. 
% Given a features matrix of $n$ examples $\bX \in \R^{n\times d}$ with $d$-dimensional features and their corresponding labels  matrix $\bY \in \R^{n\times l}$, the aim of multi-label classification problem is to assign (multiple) labels out of a total of $l$-labels to unseen test data points. The extreme multi-label classification problem refers to the setting when $n$, $d$, and $l$ quickly scale to large numbers. See \cite{Bhatia16} for descriptions and links of various extreme classification datasets.
% % The webpage \url{http://manikvarma.org/downloads/XC/XMLRepository.html} has descriptions and links of various extreme classification datasets. 
% We observe that both $\bX$ and $\bY$ are extremely sparse matrices, and $\bY$ has only ones for known labels and zeros otherwise. There exists many approaches. The embedding-based approaches have shown good results in particular. In an embedding based model, we learn our model by solving a certain optimization problem with constraints on the model size. In this paper, we exploit the structure of these optimization models as an optimization problem on matrix manifold, and propose fast solvers for both global and local embedding approach that scales for extreme classification problems.  

\section{Notation \label{notation}} In the following, $\mathbb{R}^n$ denotes space of $n$ dimensional vectors with real entries, $\mathbb{R}^{m \times n}$ denotes the space of $m \times n$ matrices with real entries, $\mathbb{R}_*^{m \times n}$ denotes the space of $m\times n$ real matrices of full rank, $\text{rank}(A)$ denotes the rank of a matrix $A,$ $\text{tr}(A)$ denote the trace of matrix $A,$ $\| x \|$ denotes the 2-norm, $\text{dim} \, U$ denotes the vector space dimension of $U,$ $\| A \|_{F}$ denotes the Frobenius norm of the matrix $A.$ For a symmetric matrix A, $A\succeq 0$ denotes that $A$ is positive semidefinite.

\section{Previous Work \label{previous_work}}

In the past few years, many novel methods have been formulated for the extreme classification problem. Most of these methods try to deal with the specific problems that are faced in the extreme classification scenario, such as scaling with large dimensions and issues with tail labels. These can broadly be classified into 4 categories.

\subsubsection{One-vs-all} strategies involve splitting the classification task by training multiple binary classifiers. State-of-the-art 1-vs-All approaches like {\tt DiSMEC} \cite{DiSMEC}, {\tt PPDSparse} \cite{PPD_Sparse} and {\tt Bonsai} \cite{bonsai} give a good prediction accuracy when compared to other state-of-the-art methods, but they are computationally expensive to train and they also have high memory requirement during training. {\tt Parabel} \cite{parabel} overcomes this problem by reducing the number of training points in each one vs all classifier. 
% Its formulation is such that each label’s negative training examples are those with most similar, or confusing relevant labels only.

\subsubsection{Tree-based Methods} aim to partition the feature and label space recursively, to break down the original problem into more feasible small scale sub-problems. Examples include {\tt SwiftXML} \cite{swiftxml}, {\tt PfastreXML} \cite{pfastrexml}, {\tt FastXML} \cite{fastxml}, and {\tt CRAFTML} \cite{craftml_siblini}. These methods generally have fast training and prediction time, but suffer in prediction performance.

\subsubsection{Embedding-based Methods} rely on low-dimensional embeddings to approximate the label space. 
The paper \cite{cssp_bi} selects a small subset of class labels such that it approximately spans the original label space. The subset construction is done by performing an efficient randomized sampling procedure where the sampling probability of each label reflects its importance among all labels. The label selection approach has previously been attempted in \cite{moplms_balasubramanian}. The label selection approach is based on the assumption that 
a small subset of labels can be used to recover all the output labels.
To select this subset, an optimization problem is solved. However, the size of the label subset cannot be controlled explicitly. To address this issue, the label subset selection problem is modelled as a column subset selection problem (CSSP). In \cite{leml_yu}, the problem of extreme classification with missing labels (i.e., labels which are present but not annotated) is addressed by formulating the problem in a generic empirical risk minimization (ERM) framework. The paper \cite{rembrandt_mineiro} utilizes a relationship between low dimensional label embeddings and rank constrained estimation to develop a fast label embedding algorithm.
In \cite{cplst_chen}, a label space dimension reduction (LSDR) approach is explored, which considers both the feature and the label parts. This approach minimizes an upper bound of the popular Hamming loss, and is called the conditional principal label space transformation (CPLST). In \cite{plst_tai}, the label space of a multi-label classification problem is perceived as a hypercube. It shows that algorithms such as binary relevance (BR) and compressive sensing (CS) can also be derived from the hypercube view. In this paper, the multilabel classification problem is reduced to binary classification in which the labelsets are represented by low-dimensional binary vectors.  This low-dimensional representation is based on the principle of Bloom filters, which is a space efficient data structure that was originally designed for approximate membership testing. The paper \cite{cs_hsu} provides some guidelines that can be used for selecting the compression and reconstruction functions for performing compressed sensing (CS) in the extreme multi-label setting. The application of CS to the XML problem is motivated by the observation that even though the label space may be very high dimensional, the label vector for a given sample is often sparse. This sparsity of a label vector will be referred to as the \textit{output sparsity.} The aim of this paper is to utilize the sparsity of $\mathrm{E}[y|x]$ rather than that of $y.$ In \cite{mach_medini}, the authors present Merged Average Classifiers via Hashing(MACH), which is a generic $L$-classification algorithm where memory scales only at $O(\log L).$ MACH is a count-min sketch structure and it relies on universal hashing to reduce classification with a large number of classes to a few independent classification tasks and with a small number of classes, hence leading to a embarrassingly parallel approach.
Earlier attempts on global embedding models, such as {\tt LEML} \cite{leml_yu}, failed to perform very well, because the tail labels were not captured accurately. More recent attempts, such as {\tt SLEEC} \cite{sleec_bhatia} and {\tt AnnexML} \cite{annexml}, used multiple local low-dimensional embeddings to circumvent this problem, and were able to achieve much better performance. Other recent methods such as {\tt ExMLDS-(4,1)} \cite{exml4} and {\tt DXML} \cite{dxml} also use embedding based approaches.

\subsubsection{Deep Learning Based Methods} have recently gained attention. These generally rely on Deep Learning models and have very good performance, but also are computationally very expensive. Examples include {\tt AttentionXML} \cite{attentionxml}, {\tt XML-CNN} \cite{xmlcnn}, {\tt DECAF} \cite{decaf}, {\tt X-Transformer} \cite{xtransformer}, {\tt LightXML} \cite{lightxml} etc. 

To the best of our knowledge, this is the first time that the Riemannian structure for the constrained optimization problem corresponding to extreme multi-label classification has been proposed.

\section{The Optimization Model \label{optimization_model}}

\subsection{Global Embedding Model}
We model the extreme classification problem with global embedding as follows
\begin{equation}
\min\limits_{\rank(W) = r} \norm{X W - Y}_{\rm{F}}^2.
\label{eq:classification_model}
\end{equation}
% where we treat no labels in $Y$ as negative labels. 
We call this the global embedding model of Extreme Classification (or simply the global model). The rank constraint with hyper parameter $r$ on the model $W$ enforces a cheap model to learn. 
It is well known that the problem (\ref{eq:classification_model}) has the following closed-form (argmin) solution 
%\begin{align}
$W^\star = V_X { \Sigma}_X^{-1} { M},$ 
%\label{eq:closed_form_solution}
%\end{align} 
where $U_X  {\Sigma}_X  V_X^\top = X$  is the thin singular value decomposition of $X$ and ${M}$ is the rank-$r$ truncated singular value decomposition of $U_X^\top Y .$ Note that $W^\star$ above is of rank $r$ by construction.
This closed-form solution is theoretically elegant, but since it involves computing the SVD of $X,$ it is not scalable for extreme classification problems with millions of features and labels. Hence, we are interested in designing a cheaper solver to estimate $W^*,$ by taking advantage of the matrix manifold structure.

% \subsection{A scalable approximation}
% We propose the following approximation to (\ref{eq:closed_form_solution}) motivated by scalability considerations to large datasets. We show the proposed solution computation in Matlab notations below.  
% \begin{equation}\label{eq:approximate_solution}
% \begin{array}{lll}
% \text{Step 1:}& [\bP, \bS, \bQ] = \texttt{svds}(\bY, r) \\
% \text{Step 2:} & \bL = \bX \setminus \bP ,\qquad \text{and}\\
% \text{Step 3:} & \bR = \bQ \star\bS.
% \end{array}
% \end{equation}
% Finally, our proposed approximate solution is $ \bW_{\rm approx}^\star = \bL  \bR^\top$, which is also of rank $r$.

% The main strength of the approximation solution approach (\ref{eq:approximate_solution}) is  that it is easy to implement. Also, the machine learning problem is now reduced to simple linear algebra operations like solving a sparse linear system.  Most platforms have optimized routines for solving sparse linear systems.

% For the cases when $\bX$ is numerically full rank (this is true for many real-world datasets). And, when $n \leq d$, the solution $\bW_{\rm approx}^\star$ is \emph{exactly} equal to $\bW^\star$.

%\section{Extreme classification as optimization problem on manifold}

\subsection{Riemannian Solver for Global Embedding Model}\label{sec:global_riemannian}
We present a Riemannian approach to the global model \eqref{eq:classification_model}. First, we identify the constraint set in the global model as a smooth submanifold and develop the geometric tools required for the optimization algorithms - tangent spaces, orthogonal projectors onto the tangent spaces, a Riemannian metric, retraction, vector transport and Riemannian gradient. These tools help us construct optimization algorithms on manifolds(\cite{boumal2020intromanifolds}, \cite{AbsMahSep2008}). Then we give an expression for the Riemannian gradient of the objective function in \eqref{eq:classification_model}. Later, an optimizer that solves \eqref{eq:classification_model} is presented which exploits this underlying geometry. The time complexity for the operations involved is presented.

A simple geometry can be given to constraint set in \eqref{eq:classification_model}. We can view it as a set of fixed-rank matrices as follows
\begin{equation}\label{eq:fixed-rank_matrices}
\R_r^{d\times l} = \{W\in \RR^{d\times l}\; :\; \rank(W) = r\}. 
\end{equation} 
The low-rank property of any $W\in \R_r^{d\times l}$ can be exploited to store it compactly in SVD form as $W = U\Sigma V^T$, where $U\in \RR^{d\times r}, V\in \RR^{l\times r}$ are orthogonal matrices and $\Sigma$ is diagonal of size $r \times r$. The set $\R_r^{d\times l}$ is a smooth embedded submanifold of $\R^{d\times l}$ of dimension $r(d+l-r)$ (\cite[Prop. 2.1]{vandereycken2013lowrank}, \cite[Section 7.5]{boumal2020intromanifolds}). The tangent space at $W\in \R_r^{d\times l}$, $W = U\Sigma V^T$, is given by 
\small
\begin{align*}
	T_W \R_r^{d\times l} = \{UMV^T+ &U_pV^T+ UV_p^T: M\in\R^{r\times r},\\ 
	 &U_p \in \R^{d\times r}, V_p\in\R^{l\times r}, U^TU_p = 0, V^TV_p = 0 \}.
\end{align*} 
\normalsize
Tangent spaces are linear spaces such that for $W\in\R_r^{d\times l}$ \cite[Thm 8.30]{boumal2020intromanifolds}
\[
\text{dim} \;T_W \R_r^{d\times l} = \text{dim} \; \R_r^{d\times l} = r(d+l-r).
\]
An important tool needed for developing optimization algorithms on smooth manifolds is the orthogonal projector map onto the tangent spaces. For a given $W\in \R_r^{d\times l}$ and $W = U\Sigma V^T$, any $Z\in \R^{d\times l}$, the orthogonal projector onto $T_W\R_r^{d\times l}$ is given by 
\begin{equation}\label{eq:orthogonal-projector}
	\text{Proj}_W(Z) = P_UZP_V + P_U^\perp ZP_V + P_UZP_V^\perp,
\end{equation} 
where $P_U = UU^T, P_V = VV^T, P_U^\perp = I_d - P_U, P_V^\perp = I_l - P_V$\cite[Section 7.5]{boumal2020intromanifolds}. Note that for any $Z\in \R^{d\times l},\text{Proj}_W(Z)\in T_W\R_r^{d\times l}$ and $\text{Proj}_W(Z) = Z$ for $Z\in T_W\R_r^{d\times l}$. 
% The tangent bundle of $\R_r^{d\times l}$ is defined as the disjoint union of all tangent spaces
% \begin{align*}
%     T\R_r^{d\times l} = \{(X,V):\;X\in\R_r^{d\times l}, V\in T_X\R_r^{d\times l}\}.
% \end{align*}

The Euclidean metric(inner product) on $\R^{d\times l}$ is defined for any $A,B\in \R^{d\times l}$ as $\langle A,B \rangle = \text{tr}(A^TB)$.
One can obtain a Riemannian metric on $\R_r^{d\times l}$ by restricting the Euclidean metric to the tangent spaces of $\R_r^{d\times l}$\cite[Prop 3.45]{boumal2020intromanifolds}
\begin{equation}\label{eq:riemannian-metric}
	\langle U,V \rangle_W = \text{tr}(U^TV),\quad W\in\R_r^{d\times l},\quad U,V\in T_W\R_r^{d\times l}.
\end{equation} 
A manifold with a Riemannian metric is called a Riemannian manifold. Consequently, a submanifold with a Riemannian metric is called a Riemannian submanifold. The Riemannian metric induces a norm on tangent spaces which is defined as
\begin{equation}
\|V\|_W = \sqrt{\langle V,V \rangle_W},
\end{equation}
where $V\in T_W\R_r^{d\times l}$ and $W\in \R_r^{d\times l}$.

Often smooth submanifolds embedded in Euclidean spaces are non-linear spaces. Hence, we cannot move along these spaces by simple linear combinations as we could in Euclidean spaces. We need maps called the retraction, which allow us to move along the manifold in a given direction. For the smooth manifold $\R_r^{d\times l}$, one such retraction is 
\begin{equation}
    R_W(V) = \underset{W'\in\R_r^{d\times l}}{\argmin\;}\|W+V-W'\|^2.
\end{equation}
The solution $W^*$ of the above optimization problem is the SVD of $W+V$ truncated to rank $r$. The sum $W+V$ can be expressed as the product $Q_1\Sigma Q_2^T$ where $\Sigma$ has rank at most $2r$ and $Q_1 \in \R^{d\times 2r},Q_2 \in \R^{l\times 2r}$ are orthogonal matrices. Using this structure, $W^*$ can be computed efficiently in $14(d + l)r^2 + r^3$ flops \cite[Algo. 6]{vandereycken2013lowrank}.

% The second order optimization algorithms are often very expensive as the size of the problem increases while steepest descent based algorithms have slow convergence rates. One would like to use approximate second order information computed cheaply from finite differences of first order derivatives on manifold. To compute finite differences of first order derivatives at two different points on a manifold 

We also need a method to translate local information at one point to another point on the manifold. Vector transport serves this purpose. For any two points on the manifold $W,W'\in\R_r^{d\times l}$, the vector transport is defined for $S\in T_W\R_r^{d\times l}$ as \cite[Prop 10.60]{boumal2020intromanifolds}
\begin{align}
    T_{W\rightarrow W'}(S) = \text{Proj}_{W'}(S).
\end{align}
Note that $T_{W\rightarrow W'}(S)\in T_{W'}\R^{d\times l}$ for all $S\in T_{W}\R^{d\times l}$. The last tool we define is the Riemannian gradient of a smooth map on a submanifold of Euclidean space. The Riemannian gradient is a generalization of the classical Euclidean gradient. For a Riemannian submanifold, the Riemannian gradient of a smooth scalar-valued function $f$ at $W\in \R_r^{d\times l}$ is defined as 
\begin{equation}
    \text{grad} \, f(W) = \text{Proj}_W(\nabla f(W)),
\end{equation}
where $\nabla f(W)$ is the Euclidean gradient of $f$ considered as a smooth function on Euclidean space \cite[Prop 3.51]{boumal2020intromanifolds}.

Now that we have finished with the discussion on the geometric tools of the manifold $\R_r^{d\times l}$, we focus our attention on solving the global model in \eqref{eq:classification_model}. For the objective function in \eqref{eq:classification_model} 
\begin{equation}
    f(W) = \|XW-Y\|_F^2,
\end{equation}
the Euclidean gradient and the Riemannian gradient of $f$ is given by 
\begin{equation}
    \nabla f(W) = 2X^T(XW-Y), \;\;\;\;\text{grad}f(W) = \text{Proj}_W(\nabla f(W)).
\end{equation}

For solving our optimization problem, we use {\tt RiemannianCG} presented in Algorithm \ref{alg:rcg}, which is a generalization of the non-linear conjugate gradients algorithm for manifolds (\cite[Algo. 1]{vandereycken2013lowrank}, \cite[Algo. 1]{boumal2015rtrmcextended}).
The step length in the above algorithm is determined by performing a line search along the search direction, that is, $\alpha_k = \text{arg min}_t f(X_k + tV_{k+1})$. This is solved by using the Armijo backtracking method \cite[Algo. 2]{boumal2015rtrmcextended}. 
% \begin{algorithm}
%     \centering
%     \caption{Armijo Backtracking}\label{alg:backtracking}
%     \begin{algorithmic}[1]
%         \Require $X \in \mathcal{M}$, $V \in T_XM$
%         \State $\alpha = 1/\| V\|$
%         \While{$f(R_X(\alpha V)) > f(x) + 0.5\cdot \text{D}f(x)[\alpha V]$ ?? } 
%         \State $\alpha = 0.5 \cdot \alpha$
%         \EndWhile \\
%         \Return $\alpha$
%     \end{algorithmic}
% \end{algorithm}
The computational cost of the projection $\text{Proj}_W(Z)$ is $O(dlr)$ flops, same as the cost for computing the vector transport. It takes $O((d + l)r^2 + r^3)$ flops for computing retraction. The Euclidean gradient for our objective can be computed in $O(dr^2 + ndl)$ flops. The per-iteration complexity is $O((d+l)r^2+ndl+r^3)$ which is grows cubically with the size of the problem. So the Riemannian approach to the global model is infeasible to large problems. Hence we turn to local embedding based models.

\begin{algorithm}[H]
	\centering
	\caption{{\tt RiemannianCG (RCG)}}\label{alg:rcg}
	\scriptsize
	\begin{algorithmic}[1]
		\Require Initial point $X_1 \in \mathcal{M}$, $V_1 = 0$ and tolerance $\tau > 0$ \Comment{Input parameters}
		\For{$k = 1,2,\cdots$} 
		\State $P_k = \text{grad}f(X_k)$
		\State if $\| P_k\| \le \tau$ then break \Comment{Check convergence}
% 		\State $\overline{P}_k = T_{X_{k - 1}\rightarrow X_k}(P_k)$
% 		\State $\overline{V}_k = T_{X_{k - 1}\rightarrow X_k}(V_k)$
		\State $Z_k = P_k - T_{X_{k - 1}\rightarrow X_k}(P_k)$
		\State $\beta_k = \max(0,\langle Z_k,P_k\rangle / \langle P_k,P_k\rangle)$
		\State $V_{k+1} = -P_k + \beta_k T_{X_{k - 1}\rightarrow X_k}(V_k)$ \Comment{Search direction}
		\State $\alpha_k = {\tt LINESEARCH}(X_k,V_{k+1})$ \Comment{Find the step length}
		\State $X_{k+1} = R_{X_k}(\alpha_k V_{k+1})$
		\EndFor
	\end{algorithmic}
\end{algorithm}

\subsection{Locally Sparse Embedding for Extreme Classification}
Earlier embedding based methods assumed that the label matrix has low rank. With this assumption, the effective number of labels is reduced by projecting the high dimensional label vector onto a low dimensional linear subspace. We note here that the low rank assumption can be violated in most real world datasets. 
% The main contribution of the paper \cite{sleec_bhatia} is a formulation for learning a small ensemble of local distance preserving embeddings which can accurately predict infrequently occurring (tail) labels. 
{\tt SLEEC} \cite{sleec_bhatia} tries to solve this issue by learning a small ensemble of local distance-preserving embeddings which can accurately predict infrequently occurring (tail) labels. 
In the earlier  embedding based methods, for a given data point $x_{i} \in \R^{d}$ an $l-$dimensional label vector $y_{i}$ is projected onto a lower $r-$dimensional linear subspace as $z_{i} = Uy_{i}.$ Regressors are then trained to predict $z_{i} = Vx_{i}.$ For a test point $x,$ prediction is done as $y = U^{\dagger}Vx,$ where $U^{\dagger}$ is a decompression matrix that takes the embedded label vectors back to the original label space. The embedding methods differ mainly in their choice of compression and decompression techniques. However, these global low rank approximations do not take into account the presence of tail labels. 
%In Fig \ref{fig:tail_labels}, it can be seen that almost $2 \times 10^{5}$ labels occurs %in less than 10 documents.
Hence, these labels will not be captured in a global low dimensional projection.
% \small
% \begin{figure}[ht]
%     \centering
%      \includegraphics[scale=0.8]{approximation_error.png}
%     \caption{\label{fig:approx_error}Plot showing the approximation error when different methods are applied for the Wiki-10 dataset. It can be seen that embeddings using low-rank approximation has a much higher error than embedding using the SLEEC approach. }
% \end{figure}
% \small
% \begin{figure}[h]
%     \centering
%      \includegraphics[scale=0.8]{tail_labels.png}
%    
%     \caption{\label{fig:tail_labels}Plot showing the no of documents containing the tail labels for WikiLSHTC dataset.  }
% \end{figure}
%
% \begin{figure}[ht]
% \begin{subfigure}{.5\textwidth}
%   \centering
%   % include first image
% %   Delicious-200K_Learner
%   \includegraphics[width=.9\linewidth]{images/approximation_error.png}  
% % \includegraphics[width=.9\linewidth]{images/Delicious-200K_Learner.png}  
%   \caption{\label{fig:approx_error}Plot showing the approximation error when different methods are applied for the Wiki-10 dataset \cite{sleec_bhatia}.}
%   \label{fig:sub-first}
% \end{subfigure}
% \begin{subfigure}{.5\textwidth}
%   \centering
%   % include second image
%  \usepackage{}  \includegraphics[width=.9\linewidth]{images/tail_labels.png}  
%   \caption{\label{fig:tail_labels}Plot showing the no of documents containing the tail labels for WikiLSHTC dataset \cite{sleec_bhatia}.}
%   \label{fig:sub-second}
% \end{subfigure}
% %\caption{Put your caption here}
% \label{fig:fig}
% \end{figure}
%
{\tt SLEEC} \cite{sleec_bhatia} instead of projecting globally onto a low-rank subspace, learns embeddings $z_{i}.$ This embedding captures label correlations in a non-linear fashion by preserving the pairwise distances between only the closest (but not all) label vectors, i.e. $\tilde{d}(z_i,z_j) = \tilde{d}(y_i,y_j)$ only if $i \in \text{kNN}(j),$ where $\tilde{d}$ is a distance metric. 
Thus, if one of the label vectors has a tail label, it will not be removed by approximation (which is the case for low rank approximation). Subsequently, the regressors $V$ are trained to predict $z_{i} = Vx_{i}.$ Rather than using a decompression matrix, during prediction, {\tt SLEEC} uses a kNN classifier in the embedding space, using the fact that the nearest neighbours have been preserved during training. Thus, for a new point $x$, the predicted label vector is obtained using $y = \sum_{i\; :\; Vx_i \in \text{ kNN}(Vx)}y_i.$ For speedup, {\tt SLEEC} clusters the training data into $C$ clusters, learns a separate embedding per cluster and performs kNN only within the test point's cluster. Clustering can be unstable in large dimensions, hence, {\tt SLEEC} learns a small ensemble, where each individual learner is generated by a different random clustering.
% Due to the presence of tail labels, the label matrix $Y$ cannot be well approximated using a low-dimensional linear subspace. Thus it is modeled using a low-dimensional non-linear manifold. That  is,  instead  of  preserving  distances of a given label vector to all the training points, SLEEC attempts to preserve the distance to only a few nearest neighbors. That is, it
{\tt SLEEC} finds an $r-$dimensional matrix $Z = [z_1,z_2, \ldots,z_n]^T \in \R^{n \times r},$ which minimizes the following objective
\begin{align}
    \underset{Z \in \R^{n\times r}}{\min}  ||P_{\Omega}(YY^T) - P_{\Omega}(ZZ^T)||_{F}^{2},
    \label{eqn:SVP}
\end{align}
where  the index set $\Omega$ denotes the set of neighbors that we would like to preserve, i.e., $(i,j) \in \Omega$ iff $j \in N_i,$ where $N_i$ denotes the set of nearest neighbors of $i$ and is selected as $N_i = \underset{S,|S| \leq \alpha.n}{\argmax}\sum_{j \in S}(y_{i}^{T}y_j),$ the set of $\alpha.n$ points with the largest inner product with $y_i$. This prevents the non-neighboring points from entering the neighborhood of any given point. Here $P_\Omega : \R^{n \times n} \rightarrow \R^{n \times n}$ is defined as
\begin{align}
    (P_{\Omega}(YY^T))_{ij} = 
    \begin{cases}
    \langle y_i,y_j \rangle, & \text{if } (i,j) \in \Omega\\
    0, & \text{otherwise}.
    \end{cases} \label{eqn:pomega}
\end{align}
Let $M = ZZ^T.$ Then \eqref{eqn:SVP} can be modeled as a low-rank matrix completion problem as shown below
\begin{equation}
    \underset{M \succeq 0, \: \text{rank}(M) \le r}{\min}  ||P_{\Omega}(YY^T) - P_{\Omega}(M)||_{F}^{2}.
    \label{eqn:SVP_form}
\end{equation}
Equation \eqref{eqn:SVP_form} can be solved by the {\tt SVP} method \cite{svp_jain}, which guarantees convergence to a local minima. Here {\tt SVP} is a projected gradient descent method, and the projection is done onto the set of low-rank matrices. The update for $t-$th step for {\tt SVP} is given by
\begin{align*}
    M_{t+1} = \tilde{P}_{r}(M_t + \eta P_{\Omega}(YY^T - M_t)), \label{eqn:SVP_update}
\end{align*}
where $M_t$ is the  iterate for the $t$-th step, $\eta > 0$ is the step size, and $\tilde{P}_{r}(M)$ is the projection of matrix $M$ onto the set of rank-$r$ positive semi-definite (PSD) matrices. Even though the set of rank-$r$ PSD matrices is non-convex, projection onto this set can be done using the eigenvalue decomposition of $M.$ Thus, we have
\begin{equation}\label{eqn:eig_M}
    \begin{split}
        M = U_M\Lambda_MU_{M}^{T} \quad 
        \Rightarrow ZZ^T = U_{M}\Lambda_{M}^{\frac{1}{2}}\Lambda_{M}^{\frac{1}{2}}U_{M}^{T} \quad 
        \Rightarrow Z = U_{M}\Lambda_{M}^{\frac{1}{2}}
    \end{split}
\end{equation}
Let $r_{M}^{+}$ denote the number of positive eigenvalues of $M.$ Let $\hat{r} = \min(r,r_{M}^{+}).$ Then taking the top-$\hat{r}$ eigenvalues of $M,$ we have $Z = U_{M}(1:\hat{r})\Lambda_{M}^{\frac{1}{2}}({1:\hat{r}}).$ Once the embeddings $Z = [z_1,z_2, \ldots z_n]^T$ are obtained, the regressor $V$ such that $Z = VX,$ is obtained by optimizing the following objective: 
\begin{align}
    \underset{V \in \R^{r \times d}}{\min} ||Z - VX||_{F}^{2} + \lambda||V||_{F}^{2} + \mu||VX||_{1}
    \label{eqn:ADMM}
\end{align}
The $L1-$term in \eqref{eqn:ADMM} makes the problem non-smooth, and it involves both $V$ and $X$. Hence, it is solved using {\tt ADMM} \cite{admm_sprechmann}, where the variables $Q$ and $V$ are alternately optimized (as shown in \cite[Algorithm 4]{sleec_bhatia}). 
% The training algorithm is shown in Algorithm \ref{alg:sleec_train} and the test algorithm is shown in Algorithm \ref{alg:sleec_test}. 
By exploiting the special structure of the matrices involved, the eigen decomposition required for {\tt SVP} can be computed in $O(r(nr+n\bar{n}))$ \cite[Section 2.1]{sleec_bhatia}, where $\bar{n}=|\Omega|/n^2$ is the average number of neighbours. Therefore, the per-iteration complexity of {\tt SVP} is $O(nr^2+nr\bar{n})$.

\subsection{Riemannian Extreme Multi-Label Classifier ({\tt RXML})} \label{section:local_model_R}
We propose a novel formulation to extreme classification problem by solving optimization problem (\ref{eqn:SVP_form}) as an optimization problem on a manifold.
The optimization problem (\ref{eqn:SVP_form}) was reformulated by modifying the constraint set.
% We can formulate the {\tt SLEEC} optimization problem (\ref{eqn:SVP_form}) as an optimization problem on a manifold by slightly modifying the constraint set. 
The set $S_+(n,r) \equiv \{M \in \R^{n \times n}: M \succeq 0,\text{rank}(M) = r\}$ is a quotient manifold \cite{MasAbs2020}. We modify the rank constraint in (\ref{eqn:SVP_form}) to $\text{rank}(M) = r$. The Riemannian version of optimization problem in (\ref{eqn:SVP_form}) is

% We can formulate the {\tt SLEEC} optimization problem (\ref{eqn:SVP_form}) as an optimization problem on a manifold by slightly modifying the constraint set. Instead of enforcing $\text{rank}(M) \le r$, we instead say $\text{rank}(M) = r$. Then, if $S_+(n,r) = \{M \in \R^{n \times n}: M \succeq 0,\text{rank}(M) = r\}$, we get the following problem 
\begin{align}
    \underset{M \in S_+(n,r)}{\min}  ||P_{\Omega}(YY^T) - P_{\Omega}(M)||_{F}^{2} \equiv \Phi(M)
    \label{eqn:Riemannian_SLEEC}
\end{align}

Any point $M \in  S_+(n,r)$ can be represented as $M = ZZ^T$, where $Z \in \R_*^{n \times r}$ is a full rank matrix. Then, $M$ is invariant under the map $Z \rightarrow ZO$, where $O$ is an orthogonal matrix of size $r$. Hence, we can consider the set $S_+(n,r)$ to be equivalent to the set $\R_*^{n \times r}/\mathcal{O}(r)$, where $\mathcal{O}(r)$ is the set of all orthogonal matrices of size $r$. The dimension of this manifold is $nr - r(r-1)/2$ \cite{MasAbs2020}.
Identifying $M$ with $Z$, we can define the objective function in \eqref{eqn:Riemannian_SLEEC} in terms of $Z$ as $f(Z) = \Phi(ZZ^T)$. It can seen that $f$ is invariant under orthogonal transformation.

The set $\R_*^{n\times r}$ is an open subset of $\R^{n\times r}$.
Thus, $\R_*^{n\times r}$ is an open submanifold of $\R^{n\times r}$ and $T_{Z}\R_*^{n\times r} = \R^{n\times r}$ \cite[Thm. 3.8]{boumal2020intromanifolds}. The decomposition of $T_{Z}\R_*^{n\times r}$ as horizontal space $$H_{Z}S_+(n,r) = \{U\in \R^{n\times r}:U^T{Z} = {Z}^TU\},$$ and its orthogonal complement can be used to uniquely represent the tangent space of $S_+(n,r)$. The horizontal space projector for any $H\in \R^{n\times r}$ is defined at ${Z}$ as
\begin{equation}
    \text{Proj}_Z^h(H) = H - ZE, 
\end{equation}
where $E$ is the solution of the Sylvester equation $E Z^TZ + Z^TZ E = Z^TH - H^TZ$ \cite[Thm. 9]{JouBacAbsSep2010}.
Similar to the case of fixed rank embedded submanifolds discussed earlier, the Riemannian metric is the trace inner product restricted to the tangent spaces. The Riemannian gradient for a smooth real valued function $\Phi$ on $S_+(n,r)$ is $\text{grad }\Phi(ZZ^T) = \nabla f(Z)$. The vector transport for any $Z_1,Z_2\in \R_*^{n\times r}$ and $U_1\in H_{Z_1}S_+(n,r)$ is
\begin{equation}
    T_{Z_1\rightarrow Z_2}(U_1) = \text{Proj}_{Z_2}^h(U_1).
\end{equation}
We use the retraction $R_Z(tU) = Z+tU$. The retraction is defined for those values of $t$ in the neighborhood of 0 for which $Z+tU$ is full-rank.

The objective function in \eqref{eqn:Riemannian_SLEEC} written using the Hadamard product notation is $f(Z) = ||P_{\Omega}\odot(YY^T) - P_{\Omega}\odot(ZZ^T)||_{F}^{2}$, where $P_{\Omega}$ is an $n \times n$ matrix. The Euclidean gradient for this objective function can be computed as
\begin{equation}
    \nabla f(Z) = 2(P_{\Omega}\odot (ZZ^T - YY^T) + P_{\Omega}^T \odot (ZZ^T - YY^T))Z
\end{equation}
Similar to the case of fixed rank embedded submanifolds discussed earlier, we use {\tt RiemannianCG} to solve the optimization problem \eqref{eqn:Riemannian_SLEEC}, the main difference being, we solve \eqref{eqn:Riemannian_SLEEC} only for each cluster and not for the entire dataset, thus ensuring the scalability of the approach. 
We use the sparse local embedding as introduced in previous section, where the subproblem \eqref{eqn:SVP_form} is solved by {\tt RiemannianCG}.
% We use the same testing strategy as {\tt SLEEC}. 
The training algorithm is shown in Algorithm \ref{alg:rxml_train} and the test algorithm is shown in Algorithm \ref{alg:rxml_test}. 

\hspace*{-0.8cm}
\begin{minipage}{\linewidth}
\begin{minipage}{0.53\textwidth}
	\begin{algorithm}[H]
		\centering
		\caption{{\tt RXML}: Train Algorithm}\label{alg:rxml_train}
		\scriptsize
		\begin{algorithmic}[1]
			\Require \\ 
			$\mathcal{D}, r, \bar{n}, C, \lambda, \mu, \rho$ (Input parameters)
			\State \text{Partition} $X$ \text{into} $Q^1,\ldots,Q^C$ using $k$-means
			\For{$\text{each partition } Q^j$}
			\State Form $\Omega^j$ using $\bar{n}$ nearest neighbors of each label vector $y_i\in Q^j$
			\State $Z^j \gets {\tt RCG}(P_{\Omega^j}, Y^jY^{j^T},r)$
			\State $V^j \gets {\tt ADMM}(X^j, Z^j, \lambda, \mu, \rho)$ 
			\State $Z^j = V^jX^j$
			\EndFor
			\State \textbf{Return:} $\{(Q^1,V^1,Z^1),\ldots,(Q^C,V^C,Z^C)\}$
		\end{algorithmic}
	\end{algorithm}
\end{minipage}
\hfill
\begin{minipage}{0.42\textwidth}
	\begin{algorithm}[H]
		\centering
		\caption{{\tt RXML}: Test Algorithm}\label{alg:rxml_test}
		\scriptsize
		\begin{algorithmic}[1]
			\Require $x,\bar{n},p$ 
			\State $Q^\tau:$ partition closest to $x$
			\State $z\gets V^\tau x$
			\State $\mathcal{N}_z\gets \bar{n}$ nearest neighbors of $z\in Z^\tau$
			\State $P_x\gets$ empirical label dist. for point $\in \mathcal{N}_z$
			\State $y_{\text{pred}} \gets {\tt Top_p}(P_x)$
		\end{algorithmic}
	\end{algorithm}
\end{minipage}  
\end{minipage}

\subsubsection{Convergence Analysis}
Convergence of {\tt RiemannianCG} follows from the following propositions. The propositions can be found in \cite{vandereycken2013lowrank} for the case where the underlying manifold is the fixed rank manifold.

\begin{proposition}
Consider an infinite sequence of iterates $X_i$ generated by {\tt RCG} with the objective $f.$ Then, every accumulation point $X_*$ of $X_i$ is a critical point of $f$. The propositions are 
\end{proposition}

It should be noted that the above proposition does not guarantee the existence of accumulation points. To this end, we work with a modified objective function which is in practice numerically equivalent to $f$.

\begin{proposition}
Consider an infinite sequence of iterates $X_i$ generated by {\tt RCG} with the objective function $g(X) = f(X) + \mu^2(\|X^\dagger\|_F^2+\|X\|_F^2)$ for some $0<\mu<1$. Then, $\underset{{i\rightarrow \infty}}{lim}\text{grad }g(X_i) = 0$.
\end{proposition}

The proof for the propositions can be found in the supplementary material. The proof is an adaptation of the one present in \cite{vandereycken2013lowrank} for the manifold $S_+(n,r)$.

\subsubsection{Space and time complexity}
Let the average number of points in each cluster be $\tilde{n}$. The complexity of computing Euclidean gradient for each cluster is approximately $O(|\Omega|r^2+\tilde{n}^2r)$. The complexity of computing projection onto horizontal space and vector transport is $O(r^3+\tilde{n}r^2)$, where $O(r^3)$ is the cost for solving Sylvester equation of size $r$.
% , where $\zeta_r$ is the computational cost for solving Sylvester equation of size $r$. We can solve the Sylvester equation in $\zeta_r  = O(r^3)$ in the worst case. 
% The Riemannian conjugate gradient requires two vector transport operations, one retraction operation, and one Riemannian gradient computation. 
Then, the per-iteration time complexity of {\tt RiemannianCG} becomes $O(|\Omega|r^2+\tilde{n}^2r+r^3)$. 
{\tt RiemannianCG} has to maintain matrices of size $\tilde{n} \times \tilde{n}$, so the worst case space complexity is $O(\tilde{n}^2)$.
% Thus, if the clustering is sparse, $|\Omega|$ will be small and our algorithm performs favourably as compared to {\tt SVP}.

\section{Numerical experiments\label{numerical_experiments}}

\subsection{Experimental setup}
Experiments were done on several publicly available multi-label datasets (see Table \ref{tab:description_datasets}, \cite{Bhatia16}). 
We have compared our method {\tt RXML} with several state-of-the-art methods: 
a) \textit{One-vs-all}: {\tt PDSparse}, {\tt Bonsai}, {\tt PPDSparse}, {\tt DiSMEC};
b) \textit{Tree-based Methods}: {\tt CRAFTML}, {\tt FastXML}, {\tt PfastreXML};
c) \textit{Embedding-based Methods}: {\tt SLEEC}, {\tt AnnexML}, {\tt ExMLDS-4}, {\tt DXML}.
The metrics we used for comparison were the P$@k\; (k=1,3,5)$ and nDCG$@k\; (k=1,3,5)$ scores, training time, testing time and, model size. We have run our code and all other methods on a Linux machine with 20 physical cores of Intel(R) Xeon(R) CPU E5-2640 v4 @ 2.40GHz and 120 GB RAM. See supplementary material for the details on parameter settings for other methods. We have not run {\tt PPDSparse} and {\tt DiSMEC} as best performance for these methods are reported on 100 core machine.
We have not compared against the recent deep learning methods as they are computationally expensive and require several GPUs, yet the train time is usually greater than 10 hours for larger datasets(\cite{attentionxml}, \cite{lightxml}).

We have used {\tt MANOPT} library \cite{manopt} for {\tt MATLAB} to implement the local embedding based Riemannian solver.
% For the global model we have used {\tt fixedrankembeddedfactory} $\R_r^{d\times l}$ of {\tt MANOPT}. 
In particular, we have used {\tt conjugategradient} ({\tt RiemannianCG}) and {\tt symfixedrankYYfactory} $S_+(n,r)$ of {\tt MANOPT}. We have used the {\tt SLEECcode} available on \cite{Bhatia16} which implements Algorithm \ref{alg:rxml_test} and rest of the Algorithm \ref{alg:rxml_train}. 

\hspace*{-0.5cm}
\begin{minipage}{\textwidth}
\begin{minipage}[b]{0.37\textwidth}
\centering
% \rule{4.7cm}{3.8cm}
% \includegraphics[height=3.6cm]{images/maxit_plot.eps}
% \captionof{figure}{P@1 score v/s maxit of {\tt RCG} and {\tt SVP} on Datasets 4, 5}
\includegraphics[height=3.6cm]{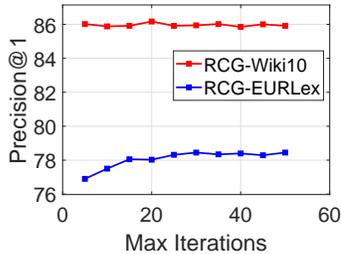}
\captionof{figure}{P@1 score v/s maxit of {\tt RCG} on Datasets with ID 4, 5.}
\label{fig:plot_maxit}
\end{minipage}
\hfill
\begin{minipage}[b]{0.63\textwidth}
\centering
\footnotesize
\begin{tabular}[H]{c|c|c|c|c|c}
\toprule
\multicolumn{1}{c}{ID} & \multicolumn{1}{c}{Dataset} & \multicolumn{1}{c}{$n$} & \multicolumn{1}{c}{$l$} & \multicolumn{1}{c}{$d$} & \multicolumn{1}{c}{$n_{test}$}\\
\midrule
	1 & {\tt Bibtex}         & 4880    & 159    & 1836   & 2515   \\
	2 & {\tt Delicious}      & 12920   & 983    & 500    & 3185   \\
	3 & {\tt Mediamill}      & 30993   & 101    & 120    & 12914  \\
	4 & {\tt EURLex-4K}      & 15539   & 3993   & 5000   & 3809   \\
	5 & {\tt Wiki10-31K}     & 14146   & 30938  & 101938 & 6616   \\
	6 & {\tt Delicious-200K} & 196606  & 205443 & 782585 & 100095 \\
	7 & {\tt Amazon-670K}    & 490449  & 670091 & 135909 & 153025 \\
	8 & {\tt AmazonCat-13K}  & 1186239 & 13330  & 203882 & 306782 \\
\bottomrule
\end{tabular}
\captionof{table}{Description of datasets.} \label{tab:description_datasets}
\end{minipage}
\end{minipage}

\subsection{Hyperparameters}
% The main hyperparameter for the global model is the approximation rank $(r)$. The rank was set by testing in the range of 10 to $L$ with suitable intervals. We found that the precision accuracy saturates by $r = 100$ on smaller datasets, which can be attributed to the inadequacy of the global model. As shown in section \ref{sec:global_riemannian}, the Riemannian approach to the global model does not scale well to larger datasets and quickly becomes infeasible. 
As we have used {\tt SLEEC} framework, the hyperparameters were set as per given in \cite[Section 3]{sleec_bhatia}. We used the hyperparameters known to perform best for {\tt SLEEC}. 
In \cite{sleec_bhatia}, the number of learners are set to be 5, 10, 20 based on the number of data points ($n$).
% For Delicious-200K and Wiki10-31K taking more than 5 learners doesn't give much increase in P$@1$ score as indicated in Fig \ref{fig:plots_learner}.
The number of clusters were chosen to be around $\left \lfloor n/6000 \right \rfloor$.
% As shown in Figure \ref{fig:plots_cluster}, taking more clusters generally did not increase performance after a threshold.
The embedding dimension is chosen to be 100 for small datasets and 50 for the large datasets. 
We have set the maximum iteration of {\tt RiemannianCG} to be 30.
% We found that our algorithm sometimes performed better even with smaller rank, as indicated in Figure \ref{fig:plots_approx_rank}.
% Fig \ref{fig:plot_cluster_iter} indicates the number of iterations taken by {\tt RCG} and {\tt SVP} to terminate on some clusters of {\tt Delicious-200K}. It was observed that {\tt SVP} terminates very early due to insufficient decrease in the objective and the error is still high ($\approx 0.1$) compared to the tolerance. In case of {\tt RCG}, the error keeps decreasing sufficiently and it seems to approach $\tau.$ Despite this fact, there is not much increase in $P@1$ score, so in order to have gain in time, we set the maximum iterations to 30 for all the large datasets. 
% The default stopping criteria for {\tt conjugategradient} is to terminate when gradient norm of the objective drops below the tolerance $\tau$, which we set to ${10}^{-3}$. But in our experimentation, we have found that ${10}^{-3}$ tolerance is rarely attained. So we added an additional stopping criteria similar to the one used in {\tt SLEEC} code \cite{Bhatia16}, where the solver terminates if the change in the objective function between two consecutive iterations is less than $0.1 \tau$.

\begin{figure}
\begin{subfigure}{.32\textwidth}
  \centering
  \includegraphics[width=\linewidth]{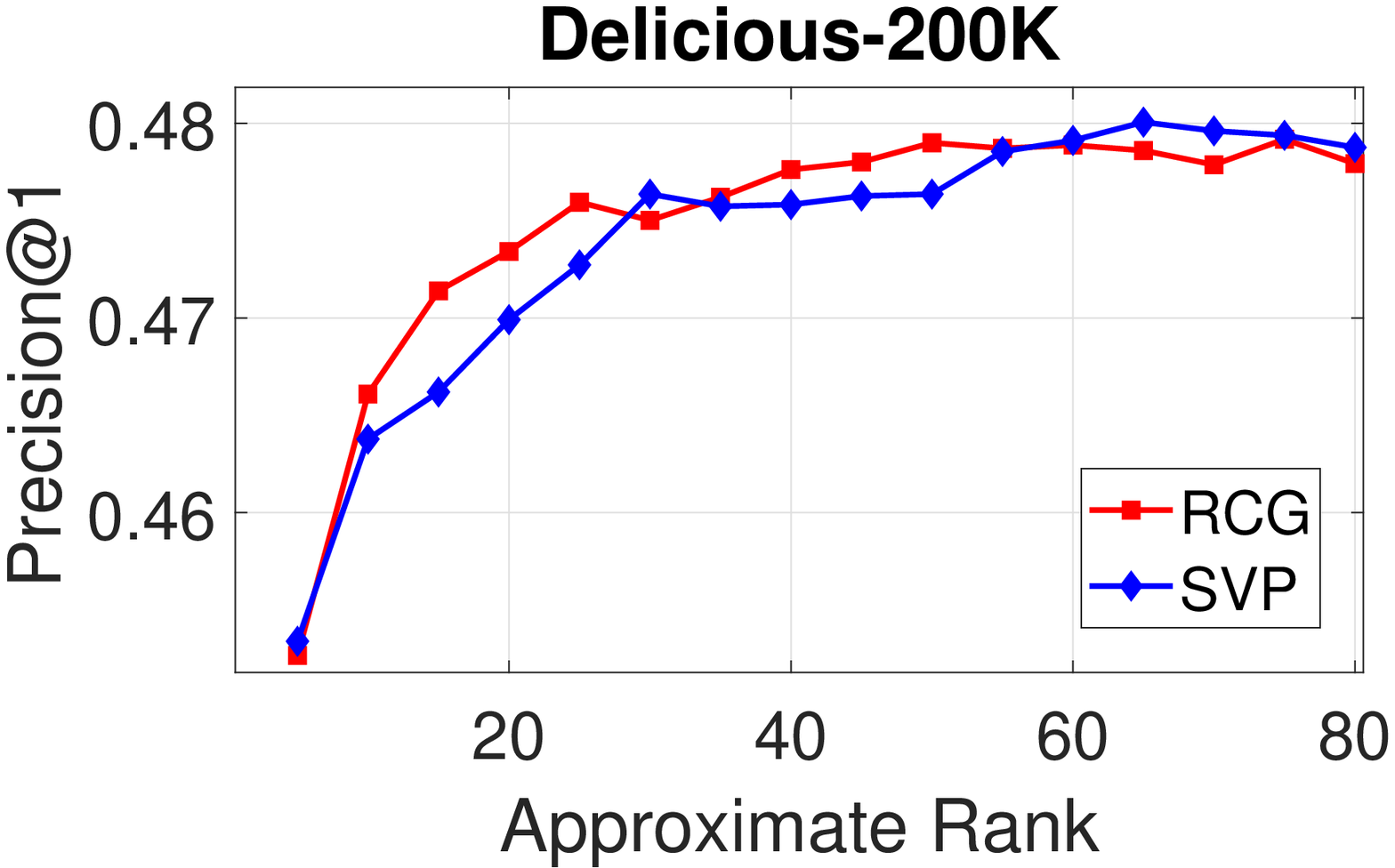}
  \caption{}
  \label{plot:rank}
\end{subfigure}
\begin{subfigure}{.32\textwidth}
  \centering
  \includegraphics[width=\linewidth]{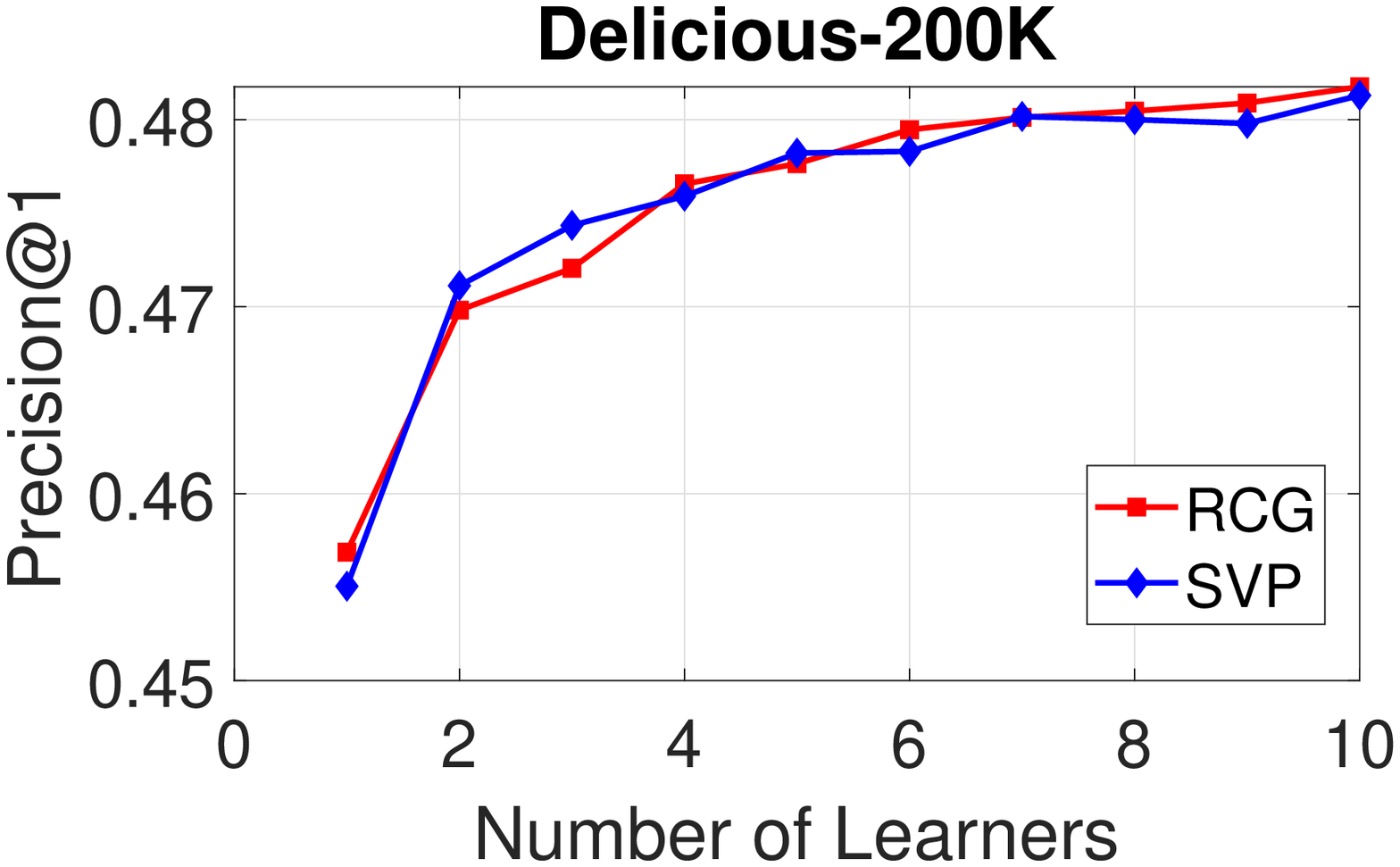}  
  \caption{}
  \label{plot:learners}
\end{subfigure}
\begin{subfigure}{.32\textwidth}
  \centering
  \includegraphics[width=\linewidth]{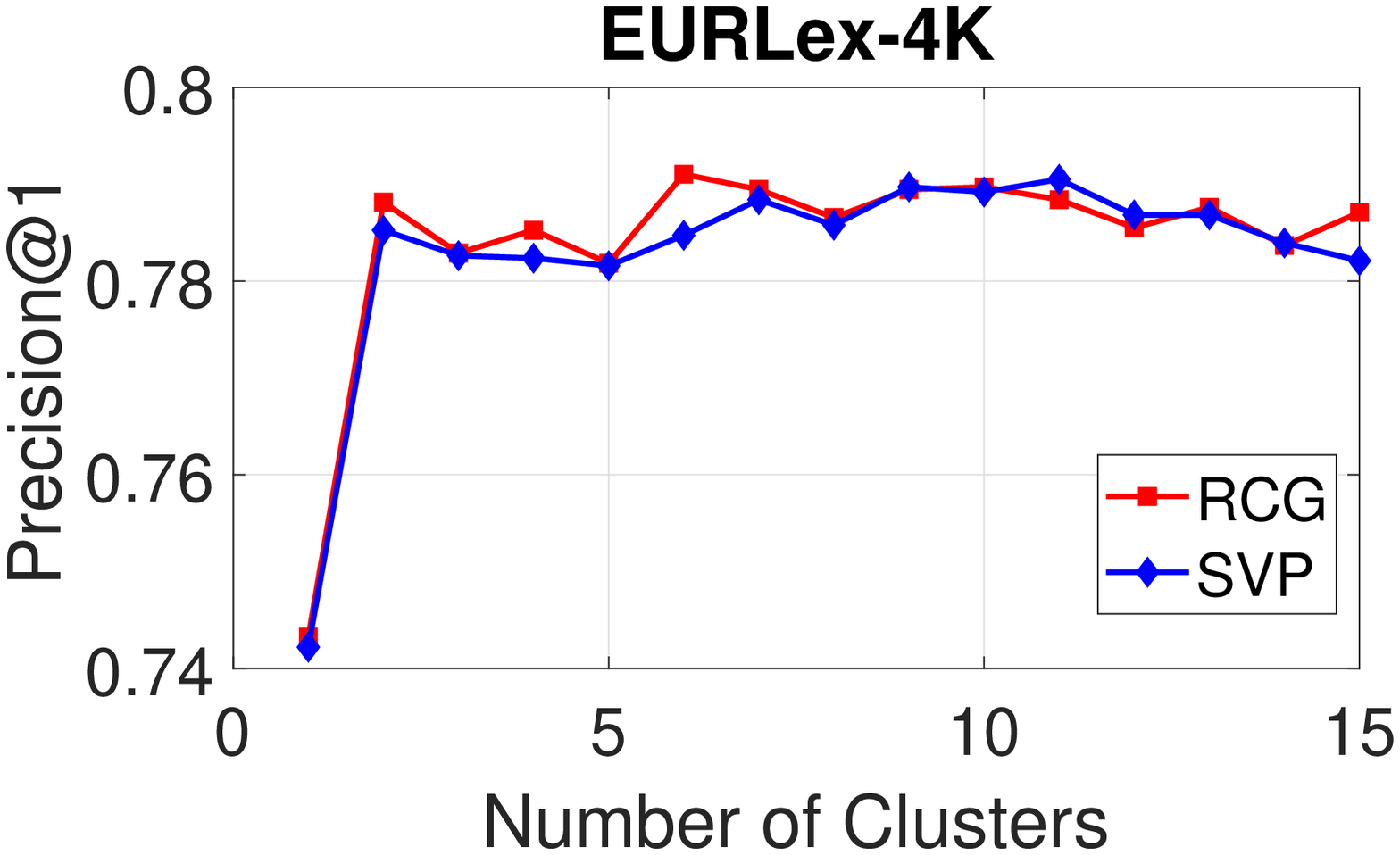}  
  \caption{}
  \label{plot:clusters}
\end{subfigure}
\caption{Variation of P@1 accuracy with embedding dimension, number of learners and number of clusters}
% The accuracy increases with the approximate rank but quickly saturates. The size of the model also increases with an increase in rank.}
\label{fig:plots_varying_params}
\end{figure}

The effect of parameters on the precision score is shown in Fig \ref{fig:plot_maxit} and \ref{fig:plots_varying_params}. Fig \ref{fig:plot_maxit} shows that the $P@1$ score is saturated after 30 iterations. Thus, setting maximum iterations to 30 gives a decent score with less train time (see Table \ref{tab:scores}, \ref{tab:times}). In our implementation we have parallelized the training over learners, so that the total train time is close to the train time of a single learner. But the model size is directly proportional to the number of learners. As indicated in Fig \ref{plot:learners}, generally 5 to 10 learners are sufficient to get good precision scores. Fig \ref{plot:rank} show that increasing embedding dimension gives a better score, but only upto a certain threshold. Further increase in score is not seen due to the limitation of the model. Choosing a lower approximation rank is crucial otherwise the model size and the per-iteration complexity of {\tt RiemannianCG} increases. Increasing the number of clusters is also crucial in reducing the per-itertaion complexity of {\tt RiemannianCG}. However as seen in Fig \ref{plot:clusters}, increasing the number of clusters does not necessarily alter the score. Hence we choose large number of clusters. 
% See supplementary for more plots showing the impact of the parameters on the score.

% Fig \ref{fig:plot_maxit} shows that after 30 maximum iterations $P@1$ score saturates for {\tt RCG} and thus justifies our choice of maxit. 

% \begin{figure}
% \begin{subfigure}[b]{.5\textwidth}
%   \centering
%   \includegraphics[width=1\linewidth]{images/time_bar_plot3.eps}  
%   \caption{}
%   \label{fig:time_bar_plot}
% \end{subfigure}
% \begin{subfigure}[b]{.5\textwidth}
%   \centering
%   \includegraphics[width=1\linewidth]{images/Delicious-200K_err_plot1.eps}
%   \caption{}
%   \label{fig:plot_cluster_iter}
% \end{subfigure}
% \caption{(a) Mean training time over learners for Riemannian and classical {\tt SLEEC} (see Table \ref{tab:description_datasets} for Dataset IDs). (b) Variation of P@1 accuracy with maximum iterations of {\tt RCG} on {\tt Wiki10-31K.}}
% \end{figure}

%% Added by pawan on 28th Feb

\begin{table}[h]
\captionsetup{font=scriptsize}
\centering
% \tiny
% \begin{adjustbox}{angle=90}
% \hspace*{-3.5cm}
% \hspace*{-0.3cm}
\scriptsize
\begin{tabular}{|c|c|c|c|c|c|c|c|c|c|c|c|c|}
  \hline
 \multicolumn{1}{|c|}{Type} & \multicolumn{5}{c|}{Embedding based} & \multicolumn{3}{c|}{Tree Based} & \multicolumn{4}{c|}{1-vs-All} \\ \hline 
            ID           & RXML  & SLEEC & AnxML & ExML4 & DXML  & CFML  & FXML  & PFXML & Bonsai& DSMC  & PPDS  & PDS   \\ \hline 
            1              &       &       &       &       &       &       &       &       &       &       &       &       \\
                       P@1 & 65.84 & 65.77 & -     & -     & 66.03 & 62.94 & 64.41 & 63.14 & -     & -     & -     & 61.29 \\
                       P@3 & 40.08 & 40.29 & -     & -     & 40.21 & 38.21 & 39.38 & 40.11 & -     & -     & -     & 35.82 \\
                       P@5 & 29.26 & 29.38 & -     & -     & 27.51 & 27.72 & 28.79 & 29.38 & -     & -     & -     & 25.74 \\ \hline 
            2              &       &       &       &       &       &       &       &       &       &       &       &       \\
                       P@1 & 67.16 & 66.91 & 64.94 & 66.17 & 68.17 & 67.88 & 68.51 & 60.31 & -     & -     & -     & 51.82 \\
                       P@3 & 61.02 & 61.12 & 59.63 & 59.87 & 62.02 & 62.00 & 63.44 & 56.08 & -     & -     & -     & 44.18 \\
                       P@5 & 56.38 & 56.40 & 55.01 & 54.85 & 57.45 & 57.32 & 59.06 & 53.41 & -     & -     & -     & 38.95 \\ \hline 
            3              &       &       &       &       &       &       &       &       &       &       &       &       \\
                       P@1 & 86.90 & 87.00 & 86.82 & 85.22 & 88.73 & 86.51 & 84.12 & 83.77 & -     & 81.86 & 86.50 & 81.86 \\
                       P@3 & 72.50 & 72.69 & 69.72 & 68.31 & 74.05 & 71.61 & 67.46 & 67.57 & -     & 62.52 & 68.40 & 62.52 \\
                       P@5 & 58.22 & 58.58 & 55.74 & 54.21 & 58.83 & 57.62 & 53.42 & 53.51 & -     & 45.11 & 53.20 & 45.11 \\ \hline 
            4              &       &       &       &       &       &       &       &       &       &       &       &       \\
                       P@1 & 79.10 & 79.36 & 80.04 & 78.41 & -     & 78.73 & 70.94 & 70.07 & 83.01 & 82.40 & 83.83 & 73.38 \\
                       P@3 & 64.37 & 64.68 & 64.26 & 61.49 & -     & 64.96 & 59.90 & 59.13 & 69.38 & 68.50 & 70.72 & 60.28 \\
                       P@5 & 52.27 & 52.94 & 52.71 & 50.30 & -     & 53.65 & 50.31 & 50.37 & 58.31 & 57.70 & 59.21 & 50.35 \\ \hline 
            5              &       &       &       &       &       &       &       &       &       &       &       &       \\
                       P@1 & 86.31 & 86.43 & 86.48 & 86.66 & 86.45 & 84.82 & 82.97 & 71.78 & 84.66 & 85.20 & -     & 77.75 \\
                       P@3 & 74.20 & 73.91 & 74.26 & 75.07 & 70.88 & 73.57 & 67.71 & 59.30 & 73.70 & 74.60 & -     & 65.46 \\
                       P@5 & 63.64 & 63.37 & 64.19 & 64.62 & 61.31 & 63.73 & 57.66 & 51.04 & 64.54 & 65.90 & -     & 55.28 \\ \hline 
            6              &       &       &       &       &       &       &       &       &       &       &       &       \\
                       P@1 & 47.86 & 47.61 & 46.64 & 47.73 & 48.13 & 47.87 & 43.09 & 17.94 & 46.69 & 45.50*& -     & -     \\
                       P@3 & 42.10 & 42.07 & 40.79 & 41.27 & 43.85 & 41.28 & 38.60 & 17.40 & 39.88 & 38.70*& -     & -     \\
                       P@5 & 39.30 & 39.35 & 37.65 & 37.93 & 39.83 & 38.01 & 36.18 & 16.99 & 36.38 & 35.50*& -     & -     \\ \hline 
            7              &       &       &       &       &       &       &       &       &       &       &       &       \\
                       P@1 & 34.79 & 34.63 & 42.08 & 38.92 & 37.67 & 37.98 & 36.60 & 36.25 & 45.58 & 44.70*& 45.32 & -     \\
                       P@3 & 31.18 & 30.90 & 36.62 & 34.30 & 33.72 & 34.05 & 32.84 & 33.56 & 40.39 & 39.70*& 40.37 & -     \\
                       P@5 & 28.60 & 28.26 & 32.77 & 30.98 & 29.86 & 31.35 & 29.93 & 31.42 & 36.60 & 36.10*& 36.92 & -     \\\hline 
            8              &       &       &       &       &       &       &       &       &       &       &       &       \\
                       P@1 & 89.87 & 89.94 & 93.54 & 92.45 & -     & 92.78*& 93.09 & 85.60 & 92.98 & 93.40*& 92.72 & 89.29 \\
                       P@3 & 75.69 & 75.82 & 78.37 & 77.43 & -     & 78.48*& 78.18 & 75.23 & 79.13 & 79.10*& 78.14 & 74.09 \\
                       P@5 & 61.28 & 61.48 & 63.31 & 62.63 & -     & 63.58*& 63.42 & 62.87 & 64.46 & 64.10*& 63.41 & 60.12 \\ \hline
% \textbf{WikiLSHTC-325K}    &       &       &       &       &       &       &       &       &      &       &       &       \\ 
%                       P@1 & 54.37 & 54.83 & 63.36 & 62.15 & -     & 56.57 & 49.75 & 56.05 & 66.6  & 64.40 & 64.08 & 61.26 \\
%                       P@3 & 33.03 & 33.42 & 40.66 & 39.58 & -     & 34.73 & 33.10 & 36.79 & 44.5  & 42.50 & 41.26 & 39.48 \\
%                       P@5 & 23.54 & 23.85 & 29.79 & 29.10 & -     & 25.03 & 24.45 & 27.09 & 33.0  & 31.50 & 30.12 & 28.79 \\ \hline 
\end{tabular}
% \end{adjustbox}
\caption{\label{tab:scores}Precision Scores of various methods. Here AnxML stands for {\tt AnnexML}, ExML4 stands for {\tt ExMLDS-4}, CFML stands for {\tt CRAFTML}, FXML stands for {\tt FastXML}, PFXML stands for {\tt PfastreXML}, DSMC stands for {\tt DiSMEC}, PPDS stands for {\tt PPDSparse}, PDS stands for {\tt PDsparse}. The missing entries indicate that the results were not reported in the source. * indicates that the particular method did not run successfully on our machine, so the reported score is taken from the source.}
\end{table}

\subsection{Results}
It is to be noted that for the methods indicated in Table \ref{tab:times}, the scores obtained during their execution are reported in Table \ref{tab:scores}. In other cases the scores are copied for either from \cite{Bhatia16} or the source paper.
\subsubsection{Comparison with Embedding-based methods}
The results for {\tt SLEEC} in Table \ref{tab:times} show orders of magnitude speed-up in train and test times compared to that reported in \cite{PPD_Sparse}. This is done by simply parallelizing the training and prediction phase over learners. {\tt RXML} is 1.5 to 2 times faster than {\tt SLEEC} in train time on large datasets with similar precision scores. In comparison with {\tt AnnexML} and {\tt ExMLDS-4}, {\tt RXML} train times are faster by 1.5 to 4 times on large dataset, and {\tt RXML} model size smaller by 3 to 10 times except for {\tt EURLex-4K}. The comparison with {\tt DXML} is ambiguous as the code is not published and it runs on GPU. Even then the train times of {\tt DXML} are reported to be much higher on {\tt Wiki10-31K}, {\tt Delicious-200K} and {\tt Amazaon-670K} with similar scores on all datasets except for {\tt Amazon-670K}.

\subsubsection{Comparison with Tree-based methods}
Again in comparison to the results reported in \cite{PPD_Sparse}, {\tt FastXML} and {\tt PfastreXML} run orders of magnitude faster with similar scores on 20 cores.
{\tt RXML} is very competitive with the state-of-the-art Tree-based methods. On {\tt EURLex-4K} and {\tt Wiki10-31K} both {\tt FastXML} and {\tt PfastreXML} are faster in terms of train times, but {\tt RXML} scores achieves a better score. {\tt CRAFTML} is slower in train time on {\tt Delicious-200K} and the prediction times are generally 10 to 15 times slower than {\tt RXML}. The memory constraints prevents us from running {\tt CRAFTML} on {\tt AmazonCat-13K}. 

\subsubsection{Comparison with 1-vs-All methods}
{\tt PDSparse} is slow compared to {\tt RXML} in train time, but it is slightly faster in prediction. {\tt Bonsai} is slower in both train and test times compared to {\tt RXML}. Memory constraints prevented us from running {\tt PDSparse} on {\tt Delicious-200K} and {\tt Amazon-670K}, and {\tt Bonsai} on {\tt Delicious-200K}, {\tt Amazon-670K}, and {\tt AmazonCat-13K}.

\subsubsection{Comparison with Deep Learning methods}
The state-of-the-art deep learning methods such as {\tt AttentionXML}, {\tt LightXML} and {\tt X-Transformer} achieve the highest precision scores, but the train times are also significantly higher compared to {\tt RXML}. In \cite{lightxml}, it is shown that {\tt AttentionXML} takes 26 hours to train on {\tt Amazon-670K}, whereas {\tt LightXML} takes 28 hours. 

\begin{table}[h]
\captionsetup{font=scriptsize}
    \centering
\scriptsize
% \tiny
% \begin{adjustbox}{angle=90}
\hspace*{-0.4cm}
\begin{tabular}{|c|c|c|c|c|c|c|c|c|c|c|}\hline
 \multicolumn{1}{|c|}{Type} & \multicolumn{5}{c|}{Embedding based} & \multicolumn{5}{c|}{Other}  \\ \hline 
    \multicolumn{1}{|c|}{Language} & \multicolumn{2}{c|}{MATLAB \& C} & \multicolumn{2}{c|}{C++} & \multicolumn{1}{c|}{Python} & \multicolumn{1}{c|}{Java} & \multicolumn{4}{c|}{C/C++} \\ \hline 
    % \multicolumn{1}{|c|}{Machine}  & \multicolumn{4}{c|}{20 cores} & \multicolumn{1}{c|}{ 50 cores} & \multicolumn{4}{c|}{1 core}\\\hline 
      Data ID              & RXML & SLEEC& AnxML& EXML4& DXML*& CFML & FXML & PFXML & Bonsai & PDS  \\ \hline 
        4                  &      &      &      &      &      &      &      &       &        &      \\
                 Train (s) & 65.8 & 715  & 146  & 106  & -    & 28   & 18.3 & 19    &  73.9  & 313  \\
                 Test (ms) & 0.48 & 0.86 & 0.05 & 0.05 & -    & 6.56 & 0.47 & 0.13  &  0.99  & 0.21 \\
                 Size (MB) & 120  & 120  & 85.8 & 85.8 & -    & 30   & 218  & 255   &  23.5  & 25   \\ \hline 
        5                  &      &      &      &      &      &      &      &       &        &      \\
                 Train (s) & 132  & 194  & 798  & 748  & 3047 & 47   & 59.3 & 65.8  &  2784  & 3558 \\
                 Test (ms) & 0.68 & 0.66 & 0.13 & 0.16 & 2.46 & 15.1 & 0.69 & 2.06  &  24.4  & 0.47 \\
                 Size (MB) & 262  & 262  & 611  & 611  & -    & -    & 482  & 1071  &  108   & -    \\ \hline 
        6                  &      &      &      &      &      &      &      &       &        &      \\
                 Train (s) & 1091 & 1607 & 4382 & 4335 & 14537& 2753 & 977  & 1054  &        &      \\
                 Test (ms) & 0.68 & 0.79 & 0.15 & 0.16 & 14.6 & 11.1 & 0.68 & 1.91  &   MLE  & MLE  \\
                 Size (GB) & 1.84 & 1.84 & 10.7 & 10.7 & -    & 0.34 & 6.2  & 14.5  &        &      \\ \hline 
        7                  &      &      &      &      &      &      &      &       &        &      \\
                 Train (s) & 1095 & 2250 & 1629 & 1609 & 67038& 848  & 803  & 851   &        &      \\
                 Test (ms) & 0.71 & 0.75 & 0.07 & 0.09 & 23.3 & 7.64 & 1.12 & 1.65  &   MLE  & MLE  \\
                 Size (GB) & 6.05 & 6.05 & 6.70 & 6.70 & -    & 0.49 & 9.31 & 10.7  &        &      \\ \hline 
        8                  &      &      &      &      &      &      &      &       &        &      \\
                 Train (s) & 3097 & 6937 & 4289 & 4286 & -    &      & 1696 & 1719  &        & 5631 \\
                 Test (ms) & 0.62 & 0.64 & 0.07 & 0.08 & -    & MLE  & 0.44 & 0.50  &   MLE  & 0.42 \\
                 Size (GB) & 5.91 & 5.91 & 18.4 & 18.4 & -    &      & 18.3 & 18.9  &        & 0.001\\ \hline 
% \textbf{WikiLSHTC-325K}    &         &         &           &          &          &        &          &        &         &          \\
%                  Train (s) & 4532.12 & 7951.34 &           & -        & 3687     & 3662   & 19160    & 20070  &         & 94343    \\
%                  Test (ms) & 0.9406  & 0.9384  &           & -        & 0.0733   & 0.0805 & 1.02     & 1.47   &         & 3.89     \\
%                  Size (GB) & 3.337   & 3.337   &           & -        & 24.207   & 24.207 & 14       & 16     &         & 0.547    \\ \hline 
\end{tabular}
% \end{adjustbox}
\caption{Train Time, Test Time(per sample), and Model Size. MLE stands for Memory Limit Exceeded. MLE happens when the memory requirements exceeds 120GB RAM.}
\label{tab:times}
\end{table}

\section{Conclusion}
% A Riemannian solver was proposed for a global embedding and a sparse local embedding. The solver is significantly faster than the SVP solver proposed earlier. Similar approach may help in other models for extreme classification.
We presented a novel Riemannian approach to Extreme classification problem, and applied it to several large scale, real world datasets. The method makes use of the underlying geometry of the {\tt SLEEC} \cite{sleec_bhatia} model to perform the optimization faster, and uses a generalization of the conjugate gradient algorithm to matrix manifolds. The model performed comparably to other well known methods, while being significantly faster to train. The method also scales well with increasing number of labels. A similar structure may also exist in other Extreme Classification models, and can be exploited to develop other variants.


\begin{thebibliography}{8}

\bibitem{boumal2020intromanifolds}
Boumal, N.: An introduction to optimization on smooth manifolds, (2020). 
Accessed online from \url{http://www.nicolasboumal.net/book}

\bibitem{vandereycken2013lowrank}
Vandereycken, B.: Low-rank matrix completion by {Riemannian} optimization. 
In: SIAM Journal on Optimization, vol. 23, pp. 1214--1236 (2013). 
\doi{10.1137/110845768}

\bibitem{manopt}
Boumal, N. and Mishra, B. and Absil, P.-A. and Sepulchre, R.: {M}anopt, a {M}atlab Toolbox for Optimization on Manifolds. 
In: Journal of Machine Learning Research, vol. 15, pp. 1455--1459 (2014). 
\url{https://www.manopt.org}

\bibitem{boumal2015rtrmcextended}
Boumal, N. and Absil, P.-A.: Low-rank matrix completion via preconditioned optimization on the {G}rassmann manifold. 
In: Linear Algebra and its Applications, vol. 475, pp. 200--239 (2015).
\doi{10.1016/j.laa.2015.02.027}

\bibitem{AbsMahSep2008} 
Absil, P.-A., Mahony, R. and Sepulchre, R.: Optimization Algorithms on Matrix Manifolds. Princeton University Press, Princeton, NJ (2008).
Accessed online from \url{https://press.princeton.edu/absil}

\bibitem{MasAbs2020}
Massart, Estelle and Absil, P.-A.: Quotient Geometry with Simple Geodesics for the Manifold of Fixed-Rank Positive-Semidefinite Matrices. 
In: SIAM Journal on Matrix Analysis and Applications, vol. 41, pp. 171--198 (2020).
\doi{10.1137/18M1231389}

\bibitem{JouBacAbsSep2010}
Journ\'ee, M., Bach, F., Sepulchre, R.and Absil, P.-A.: Low-Rank Optimization on the Cone of Positive Semidefinite Matrices. 
In: SIAM Journal on Optimization, vol. 20, pp. 2327--2351 (2010).

\bibitem{Bhatia16}
Bhatia, K. and Dahiya, K. and Jain, H. and Mittal, A. and Prabhu, Y. and Varma, M.: The extreme classification repository: Multi-label datasets and code. (2016).
\url{http://manikvarma.org/downloads/XC/XMLRepository.html}

% \bibitem{cplst_chen}
% Yao-nan Chen and Hsuan-tien Lin, Feature-aware Label Space Dimension Reduction for Multi-label Classification. (2012). In: Advances in Neural Information Processing Systems, pp. 1529--1537.

\bibitem{cplst_chen}
Chen, Y., Lin, H.: Feature-aware Label Space Dimension Reduction for Multi-label Classification. In: Advances in Neural Information Processing Systems, pp. 1529--1537, (2012).


% \bibitem{plst_tai}
% Tai, Farbound and Lin, Hsuan-Tien, Multilabel Classification with Principal Label Space Transformation. (2012). In: Neural Computation, vol. 24, issue 9, pp. 2508-2542.

\bibitem{plst_tai}
Tai, F., Lin, H.: Multilabel Classification with Principal Label Space Transformation. In: Neural Computation, vol. 24, issue 9, pp. 2508-2542, (2012).

% \bibitem{cssp_bi}
% Bi, Wei and Kwok, James T., Efficient Multi-Label Classification with Many Labels. In: Proceedings of the 30th International Conference on International Conference on Machine Learning, pp. 405-413, (2013).

\bibitem{cssp_bi}
Bi, W., Kwok, J.-T.: Efficient Multi-Label Classification with Many Labels. In: Proceedings of the 30th International Conference on International Conference on Machine Learning, pp. 405-413, (2013).

% \bibitem{moplms_balasubramanian}
% Balasubramanian, Krishnakumar and Lebanon, Guy, The Landmark Selection Method for Multiple Output Prediction. In: Proceedings of the 29th International Coference on International Conference on Machine Learning, pp. 283-290, (2012).

\bibitem{moplms_balasubramanian}
Balasubramanian, K., Lebanon, G.: The Landmark Selection Method for Multiple Output Prediction. In: Proceedings of the 29th International Coference on International Conference on Machine Learning, pp. 283-290, (2012).

\bibitem{leml_yu}
Hsiang{-}Fu Yu, Prateek Jain, Inderjit S. Dhillon, Large-scale Multi-label Learning with Missing Labels. In: arXiv, 1307.5101, (2013).



% \bibitem{rembrandt_mineiro}
% Paul Mineiro, Nikos Karampatziakis, Fast Label Embeddings for Extremely Large Output Spaces, In: arXiv, 1412.6547, (2014).

\bibitem{rembrandt_mineiro}
Mineiro, P., Karampatziakis, N.: Fast Label Embeddings for Extremely Large Output Spaces, In: arXiv, 1412.6547, (2014).

% \bibitem{sketching_edo}
% Liberty, Edo, Simple and Deterministic Matrix Sketching, In: Proceedings of the 19th ACM SIGKDD International Conference on Knowledge Discovery and Data Mining, pp. 581-588, (2013).

% \bibitem{sketching_edo}
% Liberty, E.: Simple and Deterministic Matrix Sketching, In: Proceedings of the 19th ACM SIGKDD International Conference on Knowledge Discovery and Data Mining, pp. 581-588, (2013).

\bibitem{parabel}
Y. Prabhu, A. Kag, S. Harsola, R. Agrawal and M. Varma. Parabel: Partitioned label trees for extreme classification with application to dynamic search advertising. In Proceedings of the International World Wide Web Conference, Lyon, France, April 2018

\bibitem{DiSMEC}
Babbar, R., Sch{\"o}lkopf, B.: DiSMEC: Distributed Sparse Machines for Extreme Multi-label Classification, In: WSDM, (2016).

\bibitem{PPD_Sparse}
Yen, Ian E.H. and Huang, Xiangru and Dai, Wei and Ravikumar, Pradeep and Dhillon, Inderjit and Xing, Eric, PPDsparse: A Parallel Primal-Dual Sparse Method for Extreme Classification, In: Proceedings of the 23rd ACM SIGKDD International Conference on Knowledge Discovery and Data Mining, pp. 545–553, (2017).

\bibitem{bonsai}
Khandagale, S., H. Xiao and R. Babbar. “Bonsai: diverse and shallow trees for extreme multi-label classification.” Machine Learning (2020):

\bibitem{swiftxml}
Y. Prabhu, A. Kag, S. Gopinath, K. Dahiya, S. Harsola, R. Agrawal and M. Varma. Extreme Multi-label Learning with Label Features for Warm-start Tagging, Ranking \& Recommendation. In Proceedings of the ACM International Conference on Web Search and Data Mining, Los Angeles, United States, February 2018.

\bibitem{pfastrexml}
 H. Jain, Y. Prabhu, and M. Varma, Extreme Multi-label Loss Functions for Recommendation, Tagging, Ranking \& Other Missing Label Applications. In Proceedings of the ACM SIGKDD International Conference on Knowledge Discovery and Data Mining, pp 935-944, (2016).

\bibitem{craftml_siblini}
Siblini, W., Meyer, F., Kuntz, P.: CRAFTML, an Efficient Clustering-based Random Forest for Extreme Multi-label Learning, In: ICML, (2018).

\bibitem{fastxml}
Y. Prabhu, and M. Varma, FastXML: A Fast, Accurate and Stable Tree-classifier for eXtreme Multi-label Learning, In Proceedings of the ACM SIGKDD International Conference on Knowledge Discovery and Data Mining, (2014).

\bibitem{annexml}
Tagami, Y.: AnnexML: Approximate Nearest Neighbor Search for Extreme Multi-label Classification, In: Proceedings of the ACM SIGKDD International Conference on Knowledge Discovery and Data Mining, pp 455-464, (2017).


\bibitem{attentionxml}
You, R., Zihan Zhang, Ziye Wang, Suyang Dai, H. Mamitsuka and Shanfeng Zhu. “AttentionXML: Label Tree-based Attention-Aware Deep Model for High-Performance Extreme Multi-Label Text Classification.” NeurIPS (2019).

\bibitem{exml4}
Gupta, V., Wadbude, R., Natarajan, N., Karnick, H., Jain, P., \& Rai, P. Distributional Semantics Meets Multi-Label Learning. In Proceedings of the AAAI Conference on Artificial Intelligence, 33(01), 3747-3754, (2019) 

\bibitem{dxml}
Wenjie Zhang, Junchi Yan, Xiangfeng Wang, and Hongyuan Zha. 2018. Deep Extreme Multi-label Learning. In Proceedings of the 2018 ACM on International Conference on Multimedia Retrieval (ICMR '18). Association for Computing Machinery 100–107.

\bibitem{xmlcnn}
Liu, J., Wei-Cheng Chang, Yuexin Wu and Yiming Yang. “Deep Learning for Extreme Multi-label Text Classification.” Proceedings of the 40th International ACM SIGIR Conference on Research and Development in Information Retrieval (2017)

\bibitem{decaf}
Anshul Mittal, Kunal Dahiya, Sheshansh Agrawal, Deepak Saini, Sumeet Agarwal, Purushottam Kar, and Manik Varma. DECAF: Deep Extreme Classification with Label Features. In Proceedings of the 14th ACM International Conference on Web Search and Data Mining (WSDM '21). (2021)

\bibitem{xtransformer}
Wei-Cheng Chang, Hsiang-Fu Yu, Kai Zhong, Yiming Yang, and Inderjit S. Dhillon. 2020. Taming Pretrained Transformers for Extreme Multi-label Text Classification. In Proceedings of the ACM SIGKDD International Conference on Knowledge Discovery \& Data Mining (KDD '20).163–3171.

\bibitem{lightxml}
Jiang, Ting, De-qing Wang, Leilei Sun, Huayi Yang, Zhengyang Zhao and Fuzhen Zhuang. “LightXML: Transformer with Dynamic Negative Sampling for High-Performance Extreme Multi-label Text Classification.” (2021)

% \bibitem{xreg_prabhu} 
% Prabhu, Y., Kusupati, A., Gupta, N. and Varma, M.: Extreme Regression for Dynamic Search Advertising, In: Proceedings of the 13th International Conference on Web Search and Data Mining, (2020).

\bibitem{cs_hsu}
Hsu, D.-J. , KakadeS.-M., Langford, J., Zhang, T.: Multi-Label Prediction via Compressed Sensing, In: arXiv, 0902.1284, (2009).

\bibitem{mach_medini}
Medini, T.K.R, Huang, Q., Wang, Y., Mohan, V., Shrivastava, A.: Extreme Classification in Log Memory using Count-Min Sketch: A Case Study of Amazon Search with 50M Products, In: Advances in Neural Information Processing Systems, pp. 13265-13275, (2019).

\bibitem{sleec_bhatia}
Bhatia, K., Jain, H., Kar, P., Varma, M. and Jain, P.: Sparse Local Embeddings for Extreme Multi-label Classification, In: Advances in Neural Information Processing Systems, pp. 730-738, (2015).

\bibitem{svp_jain}
Jain, P., Meka, R., Dhillon, I.-S.: Guaranteed Rank Minimization via Singular Value Projection, In: International Conference on Neural Information Processing Systems, pp. 937-945, (2010).

\bibitem{admm_sprechmann}
Sprechmann, P., Litman, R., Yakar, T.B., Bronstein, A., Sapiro, G.: Efficient Supervised Sparse Analysis and Synthesis Operators, In: Proceedings of the 26th International Conference on Neural Information Processing Systems, pp. 908-916, (2013).

\bibitem{reml_chang}
Xu, Chang and Tao, Dacheng and Xu, Chao, Robust Extreme Multi-Label Learning, In: Proceedings of the 22nd ACM SIGKDD International Conference on Knowledge Discovery and Data Mining, pp. 1275–1284, (2016).


\end{thebibliography}
\end{document}